\documentclass{article}
\usepackage{amsmath,amsthm, amsfonts, amssymb}
\usepackage{microtype}
\usepackage{fullpage}
\usepackage[square,numbers]{natbib}
\newcommand{\Z}{\mathbb{Z}}
\newcommand{\C}{\mathbb{C}}
\newcommand{\F}{\mathbb{F}}
\newcommand{\FT}{\mathcal{F}}
\DeclareMathOperator{\supp}{\mathbf{supp}}

\newtheorem{theorem}{Theorem}

\author{Saad Quader\thanks{University of Connecticut.} \and Alexander
  Russell\thanks{University of Connecticut. This research partially
    supported by the NSF under grant 1717432.} \and Ravi
  Sundaram\thanks{Northeastern University. This research partially
    supported by NSF grants 1535929 and 1718286.}}

\title{Small-Support Uncertainty Principles on $\Z/p$ over Finite Fields}

\begin{document}

\maketitle

\section{Introduction}\label{sec:intro}
Uncertainty principles are a striking manifestation of local-global phenomena in Fourier analysis, with wide ranging applications from quantum mechanics~\cite{Griffiths:Quantum} to computational complexity theory~\citep{BST}. Over a finite abelian group $G$, the classical uncertainty principle asserts that any nonzero function $f: G \rightarrow \C$ satisfies
\begin{equation}\label{eq:principle}
  |\supp f| \cdot |\supp \hat{f}| \geq |G|\,,
\end{equation}
where $\hat{f}: \Z/p \rightarrow \C$ denotes the Fourier transform of
$f$, and $\supp g$ denotes the support of the function $g$. It is easy
to check that inequality~\eqref{eq:principle} is tight for any
characteristic function $1_H$ of a subgroup $H$ of $G$. More
specifically,~\eqref{eq:principle} is tight precisely when the
function $f$ is the scaled characteristic function of a coset of a
subgroup.

This suggests that the simple groups $\Z/p$---having no nontrivial subgroups---may play a special role in the theory and, in particular, may enjoy a stronger uncertainty principle. Indeed, a recent article of \citet{Tao} establishes a remarkable strengthening over these groups: any nonzero $f: \Z/p \rightarrow \C$ satisfies
\begin{equation}  \label{eq:strong-principle}
  |\supp f| + |\supp \hat{f}| \geq p\,.
\end{equation}
A subsequent article of \citet{Meshulam} establishes an analogous statement for the groups $\Z/p^n$ under the (necessary) condition that the function $f$ is suitably ``distant'' from any subgroup coset characteristic function.
These results indicate the possibility of a refined theory that establishes stronger inequalities so long as the function under consideration avoids certain ``defects'' described by the subgroup structure.

The analogous questions over finite fields are not as well understood. For a prime power $q$ congruent to 1 modulo $p$, the Fourier transform
can be realized over the finite field $\F_q$, as there is a principal
$p$th root of unity $\omega$ in $\F_p$. In particular, the
$\F_q$-vector space of functions $\{ f: \Z/p \rightarrow \F_q\}$ is
spanned by the $\F_q$-characters
\[
  \chi_t: z \mapsto \omega^{tz}\,,\qquad t \in \Z/p\,,
\]
and the associated change of basis is carried out by the familiar Fourier transform matrix
\[
\FT = p^{-1}  \begin{bmatrix}
    \ddots & \vdots & \ddots\\
    \cdots & \omega^{k\ell} & \cdots\\
    \ddots & \vdots & \ddots
  \end{bmatrix}\,.
\]
While the uncertainty principle~\eqref{eq:principle} is preserved in this setting, the extent to which it can be strengthened along the lines of~\eqref{eq:strong-principle} remains unclear. Motivated by a direct connection to the long-standing open question of the existence of asymptotically good cyclic codes, \citet*{Lubotzky} have carried out a comprehensive study of this question. They point out some particular settings where~\eqref{eq:strong-principle} fails, but conjecture that intermediate variants---that is, inequalities stronger than~\eqref{eq:principle} but not as strong as~\eqref{eq:strong-principle}---hold in generality.

In this note, we establish a strengthened uncertainty principle for functions $f: \Z/p \rightarrow \F_q$ with a constant support (where $p \mid q-1$). In particular we show that for any constant $S > 0$, functions $f: \Z/p \rightarrow \F_q$ for which $|\supp {f}| = S$ must satisfy $|\supp \hat{f}| = (1 - o(1))p$. (Here the $o(1)$ notation denotes a function which limits to zero in $p$.) The proof relies on an application of Szemeredi's theorem and, in fact, the celebrated improvements by Gowers~\citep{Gowers} translate into slightly stronger statements in our setting (permitting conclusions for functions possessing slowly growing support as a function of $p$). 

\begin{theorem}\label{thm:main}
  Let $p$ be a prime and $q$ be a prime power for which $p \mid
  q-1$. Then any nonzero function $f: \Z/p \rightarrow \F_q$ with
  $|\supp {f}| = m$ satisfies
  \[
    |\supp \hat{f}| \geq p - r_m(p)\,,
  \]
  where $r_m(p)$ is the size of the largest subset of $\Z/p$ avoiding arithmetic progressions of length $m$.
\end{theorem}

In preparation for the remaining discussion, we record a few facts
about arithmetic progressions modulo a prime $p$. An \emph{arithmetic
  progression} of length $m$ in $\Z/p$ is a subset of the form
$\{a + kb \mid 0 \leq k < m\}$ with cardinality $m$. For a prime $p$
and a natural number $m$ we define $r_m(p)$ to be the cardinality of
the largest subset of $\Z/p$ containing no arithmetic progression of
length $m$. A historic---and highly influential---result of~\citet{Szemeredi} 
established that $r_m(p) = o(p)$ for any constant $m$ as $p \rightarrow \infty$. 
(See the book by \citet[Theorem 10.5]{TaoVuBook}.) 
In 2001, \citet[Theorem 18.6]{Gowers} established concrete upper bounds on $r_m$:
\[
  r_m(p) \leq \frac{p}{(\log \log p)^{2^{-2^{m+9}}}}\,.
\]

\begin{proof}
  With $p$ and $q$ as indicated in the statement of the theorem,
  consider a nonzero function $f: \Z/p \rightarrow \F_q$. The Fourier
  expansion of $f$ is determined by the equality
  \[
    \begin{bmatrix}
      \hat{f}(0)\\
      \vdots\\
      \hat{f}(p-1)
    \end{bmatrix}
    =
    \underbrace{
    p^{-1}
    \begin{bmatrix}
      \ddots & \vdots & \ddots\\
      \cdots & \omega^{k\ell} & \cdots\\
      \ddots & \vdots & \ddots
    \end{bmatrix}}_{\FT}
    \cdot
    \begin{bmatrix}
      f(0)\\
      \vdots\\
      f(p-1)
    \end{bmatrix}\,,
  \]
  where $\mathcal{F}$ is the Fourier transform matrix discussed above.
  Writing
  \[
    S = \supp {f}\,, \qquad m = |S|\,, \quad\text{and}\quad
    Z = \overline{\supp \hat{f}}\,,
  \]
  we first note that every minor of $\FT$ given by the columns indexed by
  $S$ and the rows indexed by a subset $Z' \subset Z$ of size $m$ must be
  degenerate. Otherwise, considering that $f$ is nonzero precisely in
  the coordinates indexed by $S$, a non-degenerate minor would
  necessarily induce a non-zero coordinate (of $\hat{f}$) in the set
  $Z' \subset Z$ (of zeros of $\hat{f}$).

  To complete the proof, we will see that any $m \times m$ minor of
  $\FT$ indexed, say, by sets $K = \{k_1, \ldots, k_m\}$ of columns
  and $L$ of rows, must necessarily have full rank if $L$ is an
  arithmetic progression. This will complete the proof, as we conclude
  that $Z$ can contain no arithmetic progression and hence must have
  cardinality no more than $r_m(p)$, as desired. It remains to check
  that a minor indexed by the sets $K$ and $L$, as described above, is
  non-degenerate. If the rows $L$ form the arithmetic progression
  $a, a + b, \ldots, a + (m-1) b$ the determinant of the resulting
  matrix can be computed exactly:
  \begin{align*}
\det    \begin{bmatrix}
      \omega^{a k_1}  & \ldots & \omega^{a k_m}\\
      \omega^{(a+b) k_1}  & \ldots & \omega^{(a+b) k_m}\\
      \vdots &  \ddots & \vdots\\
      \omega^{(a+(m-1)b) k_1} & \ldots & \omega^{(a+(m-1)b) k_m}
    \end{bmatrix}
    &=
    \det    \underbrace{\begin{bmatrix}
      (\omega^{b k_1})^0  & \ldots & (\omega^{b k_m})^0\\
      (\omega^{b k_1})^1  & \ldots & (\omega^{b k_m})^1\\
      \vdots &  \ddots & \vdots\\
      (\omega^{b k_1})^{m-1} & \ldots & (\omega^{b k_m})^{m-1}\\
    \end{bmatrix}}_{(\dag)} \;\cdot\;\prod_{i} \omega^{ak_i}\\
    &= \left(\prod_{i < j} \omega^{bk_i} - \omega^{bk_j}\right)\cdot\left(\prod_{k \in K} \omega^{ak }\right)\,,
  \end{align*}
  where we recognize $(\dag)$ as a Vandermonde matrix. This is nonzero, considering that $b \neq 0$ and the $k_i$ are distinct.
\end{proof}

\paragraph{Remark.} Over the field $\C$, every minor of the $\Z/p$ Fourier transform matrix has full rank. This fact---originally discovered by Chebotarev in the 1930s (cf. \citet{Chebotarev})---has been given attention by many authors (cf.~\citet{Dieudonne, Frenkel, EvansIsaacs, Reshetnyak, Newman}), and---along with the same relationship applied above between degenerate minors and the uncertainty principle---is the basic algebraic fact supporting the proof of the strong uncertainty principle over $\C$ by \citet{Tao}. On the other hand, \citet{Lubotzky} provide examples of minors that are degenerate in the finite field setting. The development above provides an upper estimate for how many rows can be selected while avoiding full rank. 

\paragraph{Acknowledgements.}
We thank Shai Evra for explaining the history of the problem and for pointing out that the proof can be realized over alternate coefficient rings.


\begin{thebibliography}{14}
\providecommand{\natexlab}[1]{#1}
\providecommand{\url}[1]{\texttt{#1}}
\expandafter\ifx\csname urlstyle\endcsname\relax
  \providecommand{\doi}[1]{doi: #1}\else
  \providecommand{\doi}{doi: \begingroup \urlstyle{rm}\Url}\fi

\bibitem[Griffiths(1995)]{Griffiths:Quantum}
David Griffiths.
\newblock \emph{Introduction to Quantum Mechanics}.
\newblock Pearson Prentice Hall, 2nd edition edition, 1995.

\bibitem[Barrington et~al.(1990)Barrington, Straubing, and Th{\'e}rien]{BST}
David A~Mix Barrington, Howard Straubing, and Denis Th{\'e}rien.
\newblock Non-uniform automata over groups.
\newblock \emph{Information and Computation}, 89\penalty0 (2):\penalty0
  109--132, 1990.

\bibitem[Tao(2003)]{Tao}
Terence Tao.
\newblock An uncertainty principle for cyclic groups of prime order.
\newblock \emph{arXiv preprint math/0308286}, 2003.

\bibitem[Meshulam(2006)]{Meshulam}
Roy Meshulam.
\newblock An uncertainty inequality for finite abelian groups.
\newblock \emph{European Journal of Combinatorics}, 27\penalty0 (1):\penalty0
  63--67, 2006.

\bibitem[Evra et~al.(2017)Evra, Kowalski, and Lubotzky]{Lubotzky}
Shai Evra, Emmanuel Kowalski, and Alexander Lubotzky.
\newblock Good cyclic codes and the uncertainty principle.
\newblock \emph{arXiv preprint arXiv:1703.01080}, 2017.

\bibitem[Gowers(2001)]{Gowers}
William~T Gowers.
\newblock A new proof of {Szemer{\'e}di's} theorem.
\newblock \emph{Geometric \& Functional Analysis GAFA}, 11\penalty0
  (3):\penalty0 465--588, 2001.

\bibitem[Szemer{\'e}di(1975)]{Szemeredi}
Endre Szemer{\'e}di.
\newblock On sets of integers containing k elements in arithmetic progression.
\newblock \emph{Acta Arithmetica}, 27\penalty0 (1):\penalty0 199--245, 1975.

\bibitem[Tao and Vu(2006)]{TaoVuBook}
Terence Tao and Van~H Vu.
\newblock \emph{Additive combinatorics}, volume 105.
\newblock Cambridge University Press, 2006.

\bibitem[Stevenhagen and Lenstra(1996)]{Chebotarev}
Peter Stevenhagen and Hendrik~Willem Lenstra.
\newblock Chebotar{\"e}v and his density theorem.
\newblock \emph{The Mathematical Intelligencer}, 18\penalty0 (2):\penalty0
  26--37, 1996.

\bibitem[Dieudonn{\'e}(1970/71)]{Dieudonne}
Jean Dieudonn{\'e}.
\newblock Une propri{\'e}t{\'e} des racines de l’unit{\'e}.
\newblock In \emph{Collection of articles dedicated to Alberto Gonz{\'a}lez
  Dom{\'i}nguez on his sixty-fifth birthday}, volume~25, pages 1--3. Uni{\'o}n
  Matem{\'a}tica Argentina, 1970/71.

\bibitem[Frenkel(2003)]{Frenkel}
PE~Frenkel.
\newblock Simple proof of {Chebotarev's} theorem on roots of unity.
\newblock \emph{arXiv preprint math/0312398}, 2003.

\bibitem[Evans and Isaacs(1976)]{EvansIsaacs}
RJ~Evans and IM~Isaacs.
\newblock Generalized {Vandermonde} determinants and roots of unity of prime
  order.
\newblock \emph{Proceedings of the American Mathematical Society}, 58\penalty0
  (1):\penalty0 51--54, 1976.

\bibitem[Reshetnyak(1955)]{Reshetnyak}
Yuri~Grigorievich Reshetnyak.
\newblock New proof of a theorem of {NG Chebotarev}.
\newblock \emph{Uspekhi Matematicheskikh Nauk}, 10\penalty0 (3):\penalty0
  155--157, 1955.

\bibitem[Newman(1976)]{Newman}
Morris Newman.
\newblock On a theorem of {\v{c}ebotarev}.
\newblock \emph{Linear and Multilinear Algebra}, 3\penalty0 (4):\penalty0
  259--262, 1976.

\end{thebibliography}

\end{document}